**Lazar VELIMIROVIĆ[1], Zoran PERIĆ[2], Miomir STANKOVIĆ[3], Nikola SIMIĆ[2]**

Mathematical Institute of the Serbian Academy of Sciences and Arts, Belgrade (1), Faculty of Electronic Engineering, Niš (2),
Faculty of Occupational Safety, Niš (3)

# Design of companding quantizer for Laplacian source using the approximation of probability density function

*Abstract*. *In this paper both piecewise linear and piecewise uniform approximation of probability density function are performed. For the probability density function approximated in these ways, a compressor function is formed. On the basis of compressor function formed in this way, piecewise linear and piecewise uniform companding quantizer are designed. Design of these companding quantizer models is performed for the Laplacian source at the entrance of the quantizer. The performance estimate of the proposed companding quantizer models is done by determining the values of signal to quantization noise ratio (SQNR) and approximation error for the both of proposed models and also by their mutual comparison.*

*Streszczenie.*

x

**Introduction**

It is well known that uniform quantizers are suitable for signals that have approximately uniform distribution [1]. Given that most of the signals do not have uniform distribution, there is a need for using nonuniform quantizers. Companding technique is one of the most common ways to implement nonuniform quantizers. It is based on applying a certain compressor function on an input signal. Optimal compressor function provides maximal signal to quantization noise ratio (SQNR) for referent variance of the input signal. Commonly used compression functions are the optimal compression function and $A$ and $\mu$ compression laws. Practical implementation of these compressor functions is complicated. To simplify implementation of these functions, we propose compressor function linearisation procedure.

Unlike [2] where quantizers with variable length of codewords were analyzed, this paper discusses quantizers with a fixed length of codewords. Robustness analysis of the piecewise uniform scalar quantizer is given in [3]. In addition, a comprehensive analysis of the signal to quantization noise (SQNR) in a wide range of variance for the piecewise uniform scalar quantizer for Laplacian source is derived in [4]. Approximation of the optimal compressor function using first and second order spline functions, for an input signal with Laplacian distribution, is described in [5].

In this paper the approximation of the Laplacian probability density function using linear functions is performed [6]. For the resulting approximative functions, compressor function is determined based on which the quantizer is designed and its performance is determined. As in [5], the support range of the quantizer is divided into segments of equal size, whereas the number and size of cells within segments is different. In this way we get a value of SQNR close to the value of SQNR of the nonlinear optimal companding quantizer. Depending on the method of obtaining approximative functions, in this paper piecewise linear scalar quantizers (PLSQ) and piecewise uniform scalar quantizers (PUSQ) will be discussed.

**Companding quantizer for the Laplacian probability density function**

One way of the nonuniform quantization realization is the companding technique. Nonuniform quantization can be achieved by compressing the signal $x$ using a nonuniform compressor characteristic $c(\cdot)$, by quantizing the compressed signal $c(x)$ employing a uniform quantizer, and by expanding the quantized version of the compressed signal using a nonuniform transfer characteristic $c^{-1}(\cdot)$ that is inverse to that of the compressor. The overall structure of a nonuniform quantizer consisting of a compressor, a uniform quantizer, and expandor in cascade is called compandor [1]. The granular distortion for companding quantizer is determinated by Benett's integral [1]:

(1) $$D_g = \frac{x_{\max}^2}{3N^2} \int_{-x_{\max}}^{x_{\max}} \frac{p(x)}{[c'(x)]^2} dx.$$

The function p(x) represents the probability density function of current signal values at the entrance of the scalar quantizer. In this paper, we discuss the input signal with the Laplacian distribution, which is modelled with:

(2) $$p(x) = \frac{1}{\sqrt{2}\sigma} e^{-\frac{|x|\sqrt{2}}{\sigma}},$$

where $\sigma$ is the standard deviation of the input signal $x$. The optimal compressor function $c(x)$ that achieves the maximal SQNR for referent variance of the input signal is defined as [1]:

(3) $$c(x) = x_{\max} \frac{\int_0^x p^{1/3}(x)dx}{\int_0^{x_{\max}} p^{1/3}(x)dx}, 0 \leq x \leq x_{\max}.$$

Quality of the quantized signal is estimated using SQNR that can be determined as [1]:

(4) $$SQNR = 10\log_{10}\left(\frac{\sigma^2}{D}\right). \text{[dB]}.$$

The total distortion $D$ is equal to the sum of the granular distortion $D_g$ and the overload distortion $D_o$ [1]:

(5) $$D = D_g + D_o.$$

The granular distortion is determined by the Benett's integral (1), while the overload distortion is equal to [1]:

(6) $$D_o = 2\int_{x_{\max}}^{\infty}(x - y_{\max})^2 p(x)dx.$$

where the representational level $y_{\max}$ is determined from the centroid condition:

(7) $$y_{max} = \frac{\int_{x_{max}}^{\infty} xp(x)dx}{\int_{x_{max}}^{\infty} p(x)dx} = \frac{1}{2}(\sqrt{2} + 2x_{max}).$$

**Companding quantizer construction using approximations of the probability density function**

This section gives a detailed description of two new quantizers. The first model is PLSQ, while the second proposed model is PUSQ. Support region of both quantizer models proposed in this paper is divided into $L$ segments in both quadrants, where each of the segments is divided into specified number of cells which sizes are different from segment to segment. For the realization of the quantizer which support region is formed in this way, a set of $2L$ quantizers is used. Construction of the proposed quantizer models is based on equidistance of boundary segments on which amplitude range of quantizer is divided on unequal number of cells within these segments. Approximative probability density function, on the basis of which is formed a compressor function, which is used to design PLSQ is equal to [6]:

(8) $$p^l(x) = \begin{cases} a_1 x + b_1, & x \in [0, x_1^{seg}] \\ \vdots \\ a_i x + b_i, & x \in [x_{i-1}^{seg}, x_i^{seg}] \end{cases}, i = 2,\ldots,L,$$

where coefficients $a_i$ i $b_i$ are determined as follows:

(9) $$a_i = \frac{p(x_i^{seg}) - p(x_{i-1}^{seg})}{x_i^{seg} - x_{i-1}^{seg}}, i = 1,\ldots,L,$$

(10) $$b_i = p(x_{i-1}^{seg}) - a_i x_{i-1}^{seg}, i = 1,\ldots,L.$$

Border segments are defined with:

(11) $$x_i^{seg} = i \frac{x_{max}}{L}, i = 0,1,\ldots,L.$$

The optimal support region value of the proposed quantizer, where $x_L^{seg} = x_{max}$, is as follows [7]:

(12) $$x_{max} = \frac{3}{\sqrt{2}} \ln\left(\frac{N+1}{3}\right).$$

Compressor function formed on the basis of an approximative probability density function is equal to:

(13) $$c_i^l(x) = x_{max} \frac{\int_0^x (p_i^l(x))^{\frac{1}{3}} dx}{\int_0^{x_{max}} (p_i^l(x))^{\frac{1}{3}} dx} =$$

$$= \frac{x_{max}}{\sum_{j=1}^{L} \int_{x_{j-1}}^{x_j} \frac{x_{max}}{L} (p_j^l(x))^{\frac{1}{3}} dx} \left( \sum_{j=1}^{i-1} \int_{x_{j-1}}^{x_j} (p_j^l(x))^{\frac{1}{3}} dx + \int_{x_{i-1}}^{x} (p_i^l(x))^{\frac{1}{3}} dx \right).$$

Approximative probability density function on the basis of which is formed compressor function by which we design PLSQ is defined as bellow:

(14) $$p_i^u(x) = \frac{1}{\Delta} \int_{(i-1)\Delta}^{i\Delta} p(x)dx = \frac{1}{2\Delta} \left( e^{-\sqrt{2}\Delta(i-1)} - e^{-\sqrt{2}\Delta i} \right)$$

where is $i = 1,\ldots,L$, while $\Delta$ represents the size of the segments:

(15) $$\Delta = \frac{x_{max}}{L}.$$

For an approximate probability density function formed in this way, compressor function has the following form:

(16) $$c_i^u(x) = x_{max} \frac{\int_0^x (p_i^u(x))^{\frac{1}{3}} dx}{\int_0^{x_{max}} (p_i^u(x))^{\frac{1}{3}} dx} =$$

$$= \frac{x_{max}}{\sum_{j=1}^{L} \int_{x_{j-1}}^{x_j} \frac{x_{max}}{L} (p_j^u(x))^{\frac{1}{3}} dx} \left( \sum_{j=1}^{i-1} \int_{x_{j-1}}^{x_j} (p_j^u(x))^{\frac{1}{3}} dx + \int_{x_{i-1}}^{x} (p_i^u(x))^{\frac{1}{3}} dx \right).$$

The total number of cells in the first quadrant is:

(17) $$\sum_{i=1}^{L} \frac{N_i}{2} = \frac{N-2}{2}.$$

Number of cells per segment, $N_i$, is determined from the condition:

(18) $$N_i = \frac{N-2}{2} \frac{c_i(x_i^{seg}) - c_i(x_{i-1}^{seg})}{c_L(x_L^{seg})}, i = 1,\ldots,L,$$

where for the case of PLSQ $c_i(x_i^{seg}) = c_i^l(x_i^{seg})$, while in the case of PUSQ $c_i(x_i^{seg}) = c_i^u(x_i^{seg})$. Cell threshold $x_{i,j}$, and reproduction levels $y_{i,j}$ are defined as bellow:

(19) $$x_{i,j} = c_i^{-1}(j\kappa), i = 1, j = 1,\ldots,N_i$$

(20) $$x_{i,j} = c_i^{-1}(c_i(x) + j\kappa), i = 2,\ldots,L, j = 1,\ldots,N_i$$

(21) $$y_{i,j} = c_i^{-1}\left(\left(\frac{2j-1}{2}\right)k\right), i = 1, j = 1,\ldots,N_i$$

(22) $$y_{i,j} = c_i^{-1}\left(c_i(x) + \left(\frac{2j-1}{2}\right)k\right), i = 2,\ldots,L, j = 1,\ldots,N_i$$

where for the case of PLSQ $c_i^{-1} = (c_i^l)^{-1}$, while in the case of PUSQ $c_i^{-1} = (c_i^u)^{-1}$. The step size is equal to:

(23) $$k = \frac{2x_{max}}{N-2}.$$

Granular distortion is determined by the exact formula with:

(24) $$D_g = 2 \sum_{i=1}^{L} \sum_{j=1}^{N_i/2} \int_{x_{i,j-1}}^{x_{i,j}} (x - y_{i,j})^2 p(x) dx.$$

**The performance of designed quantizer**

The performance of the proposed PLSQ model we estimate by comparing the calculated values of the signal to quantization noise ratio and error of the approximation with calculated values that correspond to the model PUSQ. Substituting equation (8), (13) and (15) in the expression (1), we obtain equations for granular distortion for PLSQ:

(25) $$D_g^l = \frac{2}{3(N-2)^2} \left( \int_0^{x_{max}} (p_i^l(x))^{\frac{1}{3}} dx \right)^3 =$$

$$= \frac{2}{3(N-2)^2} \left( \sum_{i=1}^{L} \left( \frac{3}{4a_i} \left[ (a_i \cdot i\Delta + b_i)^{4/3} - (a_i(i-1)\Delta + b_i)^{4/3} \right] \right) \right)^3,$$

where is $i = 1,\ldots,L$. Distortion overload is determined by equation (6) where the representational level $y_{max}$ is determined by (7).

PUSQ granular distortion that is obtained by combination of equation (1), (13), (16) i (18) is equal to:

(26) $$D_g^u = \frac{2}{3(N-2)^2} \left( \int_0^{x_{max}} (p_i^u(x))^{\frac{1}{3}} dx \right)^3 =$$

$$= \frac{2}{3(N-2)^2}\left(\Delta\sum_{i=1}^{L}\left(\frac{1}{2\Delta}\left(e^{-\sqrt{2}\Delta(i-1)} - e^{-\sqrt{2}\Delta i}\right)\right)^{1/3}\right)^3, i=1,\ldots,L$$

while distortion overload, as in PLSQ model, is determined using equation (6), where the representational level $y_{max}$ is determined from centroid condition (7).

In this paper, beside SQNR the error of approximation is analyzed, in order to check the accuracy of the results obtained for the values of SQNR achieved by the proposed quantizer models. The approximation error of the PUSQ model is equal to:

(27) $$e_i^u = \int_{x_{i-1}^{seg}}^{x_i^{seg}}\left|p(x)^{1/3} - p_i^u(x)^{1/3}\right|dx, i=1,\ldots,L,$$

While the approximation error of the PLSQ model equals:

(28) $$e_i^l = \int_{x_{i-1}^{seg}}^{x_i^{seg}}\left|p(x)^{1/3} - p_i^l(x)^{1/3}\right|dx. i=1,\ldots,L.$$

Total approximation error is obtained as a sum of approximation errors by segments in absolute values. For the case of PUSQ total error of approximation is equal to:

(29) $$\delta^u = \sum_{i=1}^{L}e_i^u,$$

While for the case of PLSQ total error of approximation equals:

(30) $$\delta^l = \sum_{i=1}^{L}e_i^l.$$

Table 1 shows values of SQNR and the approximation error obtained from the equation (4) for the unit variance of a signal when the number of segments of PLSQ and PUSQ is equal to $2L = 4$ and the number of quantization levels $N = 16$ and $N = 32$.

Table 1. SQNR values and the approximation error of the proposed model quantizer

| N | SQNR$^U$ [dB] | SQNR$^L$ [dB] | $\delta^u$ | $\delta^l$ |
|---|---|---|---|---|
| 16 | 17.4444 | 16.2989 | 0.3470 | 0,2694 |
| 32 | 22.5864 | 20.5028 | 0.5404 | 0,5402 |

Based on results shown in Table 1 it can be seen that the higher signal quality is achieved for the PUSQ model than for the model PLSQ. Although performances achieved for the PUSQ model are better, its implementation is much easier than for the model PLSQ. Since it was expected that with the PLSQ model we will achieve a higher quality signal then with the PUSQ model, in this paper is also analyzed error of approximation, in order to confirm obtained results for SQNR.

Based on results shown in Table 1 it can be concluded that the error of approximation is smaller for the PLSQ model than for the PUSQ model. It is shown that both based on the approximation error and based on the SQNR value, the PUSQ model presents better solution than the PLSQ model. Therefore, based on these results it can be concluded that the model PUSQ is a very efficient solution in terms of easy realization of the quantizer and achieved signal quality expressed with SQNR.

However, the following analysis shows that the results described above can not be considered as relevant.

Table 2. SQNR values of proposed quntizer models determined by the exact formula for granular distortion

| N | SQNR$^U$ [dB] | SQNR$^L$ [dB] |
|---|---|---|
| 16 | 17.4261 | 18.0277 |
| 32 | 22.6593 | 23.5937 |

Table 2 shows values of SQNR of the proposed quantizer model for the number of segments $2L = 4$ and the number of quantization levels $N = 16$ and $N = 32$ where the granular distortion is determined by the exact formula (24):

Based on the results shown in Table 2 It can be noticed that with the proposed PLSQ model higher SQNR value is achieved than with the PUSQ model, which is in contrast with results shown in Table 1. Therefore, it can be concluded that granular distortion for companding quantizers designed in this way cannot be determined using Benett's integral because a huge error can be achieved using rough approximation on which we obtain results for SQNR as shown in Table 1 and from which wrong conclusions can be drawn.

ACKNOWLEDGEMENT
This work is partially supported by Serbian Ministry of Education and Science through Mathematical Institute of Serbian Academy of Sciences and Arts (Project III44006) and by Serbian Ministry of Education and Science (Project TR32035).

***Authors***: MSc. Lazar Velimirović, Mathematical Institute of the Serbian Academy of Sciences and Arts, Kneza Mihaila 36, 11011 Belgrade, Serbia, E-mail: velimirovic.lazar@gmail.com, prof. dr Zoran Perić, Serbia, E-mail: zoran.peric@elfak.ni.ac.rs; prof. dr Miomir Stanković, Serbia, E-mail: miomir.stankovic@gmail.com; Nikola Simić, Serbia, E-mail: simicnikola90@gmail.com;